\documentclass[titlepagea4paper,11pt]{amsart}
\usepackage[arrow,matrix]{xy}
\usepackage{amssymb,yfonts}
\usepackage{amsbsy,amsmath}
\newtheorem{thm}{Theorem}[section]
\newtheorem{lemma}[thm]{Lemma}
\newtheorem{def+lem}[thm]{Definition+Lemma}
\newtheorem{cor}[thm]{Corollary}

\newtheorem{definition}[thm]{Definition}
\newtheorem{remark}[thm]{Remark}
\newtheorem{prop}[thm]{Proposition}
\newcommand{\Ext}{\mbox{\rm Ext}}
\newcommand{\re}{\mbox{\rm Re}}
\newcommand{\im}{\mbox{\rm Im}}
\newcommand{\image}{\mbox{\rm im}}

\newcommand{\dd}{\partial\bar\partial}
\newcommand{\tr}{\mbox{tr}}
\newcommand{\seq}{\longrightarrow}

\newcommand{\imp}{\Rightarrow}

\newcommand{\rk}{\mbox{\rm rk\,}}
\newcommand{\cal}{\mathcal}
\newcommand{\bb}{\mathbb}

\newtheorem{example}[thm]{Example}

\newenvironment{defn}{\begin{definition}}{\end{definition}}
\newenvironment{ex}[1]{\begin{example}\rm{\em #1}}{\hfill$\bigtriangleup$\end{example}}

\title[Ricci-Flat Deformations of Vector Bundles]{Ricci-Flat deformations of\\ holomorphic vector bundles}

\author{Marco K\"uhnel}
\address{Marco K\"uhnel\\Otto-von-Guericke-Universit\"at Magdeburg\\FMA / IAN\\Postfach 4120\\
39016 Magdeburg\\Germany}
\email{marco.kuehnel@ovgu.de}
\subjclass{14J60}
\keywords{hermitian, vector bundle, non-filtrable, deformation, Hopf manifold, curvature}

\date{\today}

\begin{document}

\begin{abstract}In this paper we give a criterion for a deformation of a holomorphic vector bundle to be Ricci-flat. As an application we show
that on a K\"ahler manifold, every deformation of a holomorphic vector bundle can be made Ricci-flat whereas on some Hopf manifolds, the non-existence of a Ricci-flat
deformation of a vector bundle is connected to non-trivial vector bundles on ${\bb C}^n\setminus\{0\}$. In the case of surfaces with $b_1(X)\not=0$, 
'Ricci-rigid' rank 2 vector bundles are proved to be strongly obstructed with respect to possible filtrations. 
We use this to conclude that on diagonal Hopf surfaces apart from the line bundles only non-good rank $2$ bundles can be Ricci-rigid. On Inoue surfaces all
Ricci-rigid rank $2$ bundles are non-filtrable except one.
\end{abstract}

\maketitle

\section{Introduction}

Let $X$ be a compact complex manifold, $U$ a complex manifold without holomorphic vector bundles (e.g. a contractible Stein manifold), 
and ${\cal E}\seq X\times U$ a holomorphic vector bundle. (So the deformations we deal here with are neither global nor small but something in between.)
Further let $g$ be a hermitian metric on ${\cal E}$. The data of the central fibre will be denoted by ${\cal E}_0$ and $g^0$.
This setup is called a deformation of the hermitian vector bundle $({\cal E}_0,g^0)$.
As usual we denote $\Theta^0:=\dd\log\det g^0$ and $\Theta^t:=\partial_X\overline\partial_X\log\det g|_{{\cal E}_t}$ the curvatures.

\begin{defn}The deformation $({\cal E},g)$ is called Ricci-flat, if $\Theta^t=\Theta^0$ for all $t\in U$.\end{defn}

We tackle the question of existence of Ricci-flat deformations. 
In Theorem \ref{main} we show that the obstruction space for this problem
is $$H^1(X,{\cal O}_X)/i_*H^1(X,{\bb R}),$$ if $i:{\bb R}\seq{\cal O}_X$ is the natural inclusion. This group is trivial for
K\"ahler manifolds, so on K\"ahler manifolds we can extend every metric on the central fibre curvature preservingly to the
deformation. In the non-K\"ahler case we can still show that the property of a given vector bundle to admit only Ricci-flat deformations is independent of the vector bundle and hence
a property of the manifold, namely the vanishing of the obstruction space. 

The answer to the question, whether for a vector bundle any Ricci-flat deformation is trivial, depends on the bundle,
however, and will be denoted by Ricci-rigidity. This generalises the notion of rigidity to manifolds with non-vanishing first cohomology. 
We will see that this property is connected to a minimizing property of $H^1({\cal End}({\cal E}))$.

On Hopf manifolds we deepen this connection and so construct Ricci-rigid vector bundles.
We relate Ricci-rigidity to the non-triviality of the pullback to the universal cover ${\bb C}^n\setminus\{0\}$. 

In the case of a surface with $b_1(X)\not=0$ Ricci-rigidity of rank 2 bundles will be shown to obstruct possible filtrations. We use this to prove that on an Inoue surface
all Ricci-rigid rank $2$ bundles are non-filtrable except one (up to tensoring with line bundles) and on Hopf surfaces given by a diagonal automorphism the only Ricci-rigid
vector bundles of rank $>1$ are of rank $2$ and not good in the sense of deformation theory.

It should be noted that the algebraic and analytic category differ widely in this case. Whereas on ${\bb C}^2$ and ${\bb C}^2\setminus\{0\}$ there are only trivial algebraic 
vector bundles \cite{se58,ho64}, on ${\bb C}^2\setminus\{0\}$ some non-trivial holomorphic
vector bundles have been constructed before. The paper by B\u anic\u a and LePotier \cite{bp87} classified filtrable holomorphic 
vector bundles on non-algebraic surfaces. In particular, for all integers $r\ge 2,c_2\ge 0$ there exists a holomorphic
rank $r$ bundle ${\cal E}$ with $c_2({\cal E})=c_2$ on a Hopf surface. Calculations
below and in \cite{ma92} show that the pullback to ${\bb C}^2\setminus\{0\}$ of such a bundle is not trivial provided $c_2>0$. 
The constructed bundles are filtrable. Later Ballico \cite{b02} constructed
a non-filtrable rank $2$ bundle on ${\bb C}^2\setminus\{0\}$. In \cite{s66} a non-trivial line bundle on ${\bb C}^2\setminus\{0\}$ is constructed.

\section{The local data}

By the compactness of $X$ we can employ the trivializations of ${\cal E}$ in order to obtain a finite cover of open sets $U_i\subset X$ 
and isomorphisms
$$\psi_i:{\cal E}|_{U_i\times U}\seq pr_1^*{\cal E}_0|_{U_i\times U},$$
where $pr_1:X\times U\seq X$ is the projection. This yields the data
$$\theta_{ij}:=\psi_i\circ\psi_j^{-1}, g_i:=(\psi_i^{-1})^*g,$$
satisfying $g_j=\theta_{ij}^*g_i$, obviously. Note that $\theta_{ij}|_{U_{ij}\times\{0\}}=Id$. 
Similarly, we have the extended metric of the central fibre
$$<,>:=pr_1^*g^0$$
as a comparison metric. We obtain by Lax-Milgram
$G_i\in End_{C^\infty}(pr_1^*{\cal E}_0|U_i\times U)$
with the property
$$g_i(e_1,e_2)=<e_1,G_ie_2>$$
for any $C^\infty$ sections $e_1,e_2$ of $pr_1^*{\cal E}_0$ over $U_i$. Note that again $G_i|_{U_i\times\{0\}}=Id$. On $U_{ij}\times U$
we have the formula
\begin{equation}\label{gs}G_i=\theta_{ji}^*G_j\theta_{ji}.\end{equation}
Here $\theta_{ji}^*$ denotes the adjoint of $\theta_{ji}$ with respect to $<,>$.
It is easy to see that the deformation is Ricci-flat if and only if 
\begin{equation}\label{det}\partial_X\overline\partial_X\log\det G_i=0\mbox{ for all }i.\end{equation}

The inclusions ${\bb Z}\subset{\bb R}\subset{\cal O}_X$ give a commuting triangle
$$\xymatrix{H^1(X,{\bb Z})\ar[r]\ar[rd]&H^1(X,{\bb R})\ar[d]^{i_*}\\&H^1(X,{\cal O}_X)}$$
of injective maps. We view all unnamed maps as natural inclusions.

For any deformation the map
$$\eta:U\seq H^1(X,{\cal O}_X)/H^1(X,{\bb Z})=Pic^0(X)$$
given by 
$$t\mapsto [\frac 1{2\pi i}\log\det\theta_{ij}(.,t)]=\det{\cal E}_t\otimes\det{\cal E}_0^\vee$$
is a well-defined holomorphic map with $\eta(0)=0$. If we assume the deformation to be small, we may assume
that $\eta$ maps to $H^1(X,{\cal O}_X)$.

\section{The infinitesimal data}

We give the infinitesimal data for the one-dimensional case $U=({\bb C},0)$. The generalisation to $U=({\bb C}^n,0)$ is
straightforward.

We take the Taylor series to the first order
$$\theta_{ij}=Id+th_{ij}+h.o.t., G_i=Id+t\gamma_i+\overline t\gamma_i^*+h.o.t.$$
with $h_{ij}\in End({\cal E}_0|U_{ij}), \gamma_i\in End_{C^\infty}({\cal E}_0|U_i)$. The $h_{ij}$ satisfy the cocycle condition
and hence yield an element $h\in H^1({\cal End}({\cal E}_0))$. Indeed, if the deformation is trivial, $h_{ij}$ is a coboundary.
By comparison of
$$g_i(e_1,e_2)=<e_1,e_2>+t<e_1,\gamma_i^*e_2>+\overline t<e_1,\gamma_ie_2>+h.o.t.$$
and
\begin{eqnarray*}\theta_{ji}^*g_j(e_1,e_2)&=&<e_1,e_2>+t(<h_{ji}e_1,e_2>+<e_1,\gamma_j^*e_2>)+\\
& &+\overline t<e_1,(\gamma_j+h_{ji})e_2>+h.o.t.\end{eqnarray*}
we obtain
\begin{equation}\label{hg}h_{ij}=\gamma_j-\gamma_i\end{equation}
on $U_{ij}$. This tells us that the $\gamma_i$ trivialise the deformation in the $C^\infty$ sense to the first order.

Now let us consider curvature
\begin{eqnarray*}\Theta_i&:=&\partial_X\overline\partial_X\log\det g_i\\
&=&\Theta^0+\dd\log\det G_i\\
&=&\Theta^0+\dd\log(1+t\tr\gamma_i+h.o.t)\\
&=&\Theta^0+t\dd\tr\gamma_i+h.o.t.\end{eqnarray*}
So, any Ricci-flat deformation satisfies
\begin{equation}\label{trace}\dd\tr\gamma_i=0.\end{equation}

\section{Existence}

We start with an arbitrary metric $\tilde g$ on a deformation ${\cal E}$. This yields the data $\tilde\psi_i$, $\tilde g^0$ and $\tilde G_i$. By the Cartan
decomposition we can find a unique positive hermitian $A\in End_{C^\infty}({\cal E}_0)$ with respect to $\tilde g^0$ such that
$$g^0(e_1,e_2)=\tilde g^0(Ae_1,Ae_2)$$
for any local sections $e_1,e_2$ of ${\cal E}_0$. Setting $\psi_i:=A^{-1}\circ\tilde\psi_i, G_i:=A^{-1}\tilde G_iA$ (now with respect to $<,>:=pr_1^*g^0$) 
we obtain a deformation of $g^0$.   

The question of existence of a Ricci-flat deformation is more subtle. 

\subsection{Obstructions}

\begin{thm}\label{main}A deformation ${\cal E}\seq X\times U$ of a hermitian holomorphic vector bundle $({\cal E}_0\seq X,g^0)$ 
allows for a curvature preserving metric $g$ extending $g^0$ if and only if 
$$\eta(t)\in i_*H^1(X,{\bb R})/H^1(X,{\bb Z})$$
for all $t\in U$ and
the inclusion $i:{\bb R}\hookrightarrow{\cal O}_X$; moreover, this condition being satisfied, in every conformal class of metrics
$g$ deforming $g^0$ there is a Ricci-flat representative. 
\end{thm}

\begin{proof}Necessity: We shrink $U_i$ appropriately such that they are simply connected. By (\ref{det}) we find for every $t\in U$
holomorphic $h_i(t)\in{\cal O}_X(U_i)$ such that
$$\log\det G_i(t)=\re(h_i(t)).$$
So we have by (\ref{gs}) on $U_{ij}$
$$2\re(\log\det\theta_{ij}(t))=\re(h_i(t)-h_j(t)),$$
hence there are $c_{ij}(t)\in{\bb R}$ such that
$$\frac{1}{2\pi i}(2\log\det\theta_{ij}(t)-h_i(t)+h_j(t))=c_{ij}(t),$$
but this means exactly
$$\eta(t)\in i_*H^1(X,{\bb R})/H^1(X,{\bb Z}).$$

Sufficiency: Since $\theta_{ij}|{U_{ij}\times\{0\}}=Id$, by shrinking $U$ we may assume that $\theta_{ij}=\exp(2\pi ik_{ij})$ for some 
$k_{ij}\in End(pr_1^*{\cal E}_0|{U_{ij}\times U})$.
Note that $(k_{ij})$ is no cocycle, unless they commute. Nevertheless, 
$\det\theta_{ij}=\exp(2\pi i\tr k_{ij})$, so $(\tr k_{ij})\in{\cal O}(U_{ij}\times U)$ is a cocycle, defined uniquely by $\theta_{ij}$
up to an integer. In particular, $(\tr k_{ij}|U_{ij}\times\{t\})\in{\cal O}(U_{ij})$ 
is a cocycle for all $t\in U$. Since we assumed that $\eta(t)=[\tr k_{ij}(t)]\in i_*H^1(X,{\bb R})/H^1(X,{\bb Z})$, 
we can choose $K_i(t)\in{\cal O}(U_i)$ and $\phi_{ij}(t)\in{\bb R}$ such that
$$\tr k_{ij}=K_j-K_i+\phi_{ij}$$
and $K_i(0)=0, \phi_{ij}(0)=0$. 
The dependence on $t$ of $K_i$ and $\phi_{ij}$ is not holomorphic anymore, but it can be chosen to be $C^\infty$. For $H_i:=\exp(2\pi iK_i)$ we obtain
\begin{equation*}\det\theta_{ij}\exp(-2\pi i\phi_{ij})=\frac{H_j}{H_i}\end{equation*}
and hence the deformed metric has to satisfy
\begin{equation}\label{detG}\det G_j=\frac{|H_j|^2}{|H_i|^2}\det G_i.\end{equation}

Now we can take any metric $\tilde g$ deforming $g^0$ (with data $\tilde G_i$) and take a conformal change:
$$G_i:=\left(\frac{|H_i|^2}{\det\tilde G_i}\right)^{\frac 1n}\tilde G_i.$$
We obtain $\det G_i=|H_i|^2$ and hence $\dd\log\det G_i=0$, so we have a Ricci-flat deformation.

Now covering the original $U$ by neighbourhoods like above and
using the uniqueness of $\tr k_{ij}$ up to an integer, we extend the criterion to all of $U$.
\end{proof}

In the infinitesimal information we loose sufficiency of the condition:

\begin{prop}\label{maxcplx}A Ricci-flat deformation germ of hermitian vector bundles satisfies
$$\frac 1{2\pi i}\tr h\in V,$$
where $V$ is the maximal complex subspace of $i_*H^1(X,{\bb R})$.
%for the inclusion $i:{\bb R}\seq{\cal O}_X$. 
\end{prop}

\begin{proof}We shrink $U_i$ appropriately such that they are simply connected. 
By (\ref{trace}) we find holomorphic
$f_i, g_i\in{\cal O}(U_i)$ such that $\tr\gamma_{i}=f_i+\overline g_i.$ So (\ref{hg}) tells us that
$$\tr h_{ij}+f_i-f_j=\overline g_j-\overline g_i.$$
Since the right hand side is antiholomorphic and the left hand side is holomorphic, we find constants $c_{ij}\in{\bb C}$
such that
$$c_{ij}=\tr h_{ij}+f_i-f_j=\overline g_j-\overline g_i.$$
So we obtain in $i_*H^1(X,{\bb C})$
$$\frac 1{2\pi i}\tr h=\frac 1{2\pi i}\tr c=\frac 1{2\pi i}(\overline c+2i\im c)=\frac{\im c}{\pi}\in i_*H^1(X,{\bb R}).$$

If we look at a base transformation $\tau:(T,0)\seq (T,0)$ and denote the 
objects corresponding to the deformation $\tau^*{\cal E}$ of ${\cal E}_0$ by a superscript $\tau$, it is straightforward that
$$h_{ij}^\tau=\tau'(0)h_{ij}.$$
This proves the statement.
\end{proof}

\section{Stable Curvature and Ricci-Rigidity}

We see that Ricci-flatness of a deformation does not depend on the initial metric on the central
bundle. There are two natural properties connected to Ricci-flat deformations:

\begin{defn}The vector bundle ${\cal E}_0$ is said to have {\em stable curvature}, if all deformations of ${\cal E}_0$ over a base $U$
without holomorphic vector bundles are Ricci-flat.
${\cal E}_0$ is {\em Ricci-rigid}, if every small Ricci-flat deformation of ${\cal E}_0$ is trivial.
\end{defn}

Our first goal is to realise that stable curvature is not a property of a vector bundle, but of the underlying manifold:

\begin{prop}\label{stcur}Let $X$ be a compact complex manifold. Then the properties
\begin{enumerate}
\item There exists a vector bundle ${\cal E}_0$ with stable curvature,
\item All vector bundles have stable curvature,
\item $i_*:H^1(X,{\bb R})\seq H^1(X,{\cal O}_X)$ is an ${\bb R}$-isomorphism
\end{enumerate}
are equivalent.
\end{prop}

\begin{proof}Since $(ii)\imp(i)$ is obvious and $(iii)\imp(ii)$ is Theorem \ref{main}, it remains to show $(i)\imp(iii)$. So 
we assume that $i_*$ is not surjective. In particular, $\dim Pic^0(X)>0$, hence ${\cal O}_X$ can be deformed non-trivially.
So, for any deformation ${\cal L}_t$ of ${\cal L}_0={\cal O}_X$ with ${\cal L}_{t}\not\in i_*H^1(X,{\bb R})/H^1(X,{\bb Z})$
for all real $t$ and arbitrary vector bundle ${\cal E}_0$, the deformation ${\cal E}_0\otimes{\cal L}_t$ satisfies
$\eta(t)={\cal L}_t^n\not\in i_*H^1(X,{\bb R})/H^1(X,{\bb Z})$ for all real $t$. Hence we have shown that no vector bundle has
stable curvature.
\end{proof}

Due to this result we can define

\begin{defn}A compact complex manifold $X$ is called a {\em stable curvature manifold}, if one (and hence all) of the properties in Proposition \ref{stcur} are satisfied.
\end{defn}

Examples of stable curvature manifolds are all compact K\"ahler manifolds. Counterexamples
are all compact manifolds with odd $b_1(X)$, e.g. Hopf surfaces.

Also Ricci-rigidity has geometric implications.

\begin{prop}\label{v=0}If $X$ allows for a Ricci-rigid vector bundle, then all germs of holomorphic maps $f:({\bb C},0)\seq Pic^0(X)$ with image
in $$i_*H^1(X,{\bb R})/H^1(X,{\bb Z})$$ are constant; in particular, Ricci-flat deformations are exactly the ones with preserved determinant; and $V=\{0\}$.
\end{prop}

\begin{proof}
% The claimed property is equivalent to all germs of holomorphic maps $f:({\bb C},0)\seq Pic^0(X)$ with image
% in $$i_*H^1(X,{\bb R})/H^1(X,{\bb Z})$$ being constant. In order to see this,
Let $f$ be such a map with $f(0)=0$ and ${\cal E}_0$ a Ricci-rigid vector bundle. 
Then, of course ${\cal E}_0\otimes f(t)$ is a Ricci-flat deformation, hence ${\cal E}_0\otimes f(t)\cong{\cal E}_0$. Taking determinants
shows $f(t)^n={\cal O}_X$, but this means that $f(t)$ is constant, so $f\equiv 0$. 
\end{proof}

Note, that for any vector bundle the sequence 
$$0\seq{\cal O}_X\stackrel{\cdot Id}{\seq}{\cal End}({\cal E}_0)\stackrel{\psi}{\seq}ad({\cal E}_0)\seq 0$$
splits via the trace map: 
$${\cal End}({\cal E}_0)\cong{\cal O}_X\oplus ad({\cal E}_0), \phi\mapsto\left(\frac {\tr\phi}{\rk{\cal E}_0},\psi(\phi)\right)$$
is an isomorphism.

In particular, $tr:H^1({\cal End}({\cal E}_0))\seq H^1({\cal O}_X)$ is surjective. 
%Let $V$ be the maximal complex subspace of $i_*H^1(X,{\bb R})\subset H^1(X,{\cal O}_X)$.
If ${\cal E}_0$ is good, i.e. $H^2(ad({\cal E}_0))=0$, then we can integrate any holomorphic $h:U\seq V$ on a small open set $U\subset{\bb C}$ to
a Ricci-flat deformation of ${\cal E}_0$. 

The classical correspondence that a vector bundle ${\cal E}_0$ over $X$ with $H^2(ad({\cal E}_0))=0$ 
is rigid if and only if 
$H^1(X,{\cal End}({\cal E}_0))=0$ has a nice analogue here. 
Our notion of Ricci-rigidity measures in the infinitesimal setup that at ${\cal E}_0$ the minimal value is attained.

\begin{prop}\label{rrcoh}If a good vector bundle ${\cal E}_0$ on $X$ is Ricci-rigid, then $H^1(ad({\cal E}_0))=0$.
\end{prop}

\begin{proof}Let ${\cal E}_0$ be Ricci-rigid and ${\cal E}_t$ be any small deformation of ${\cal E}_0$ and define $\tilde{\cal L}_t:=\det{\cal E}_t^\vee\otimes
\det{\cal E}_0$. Note that $\tilde{\cal L}_t\in Pic^0(X)$.
If we shrink the base of the deformation appropriately we may assume the existence of a holomorphic family ${\cal L}_t$
of line bundles with ${\cal L}_t^n=\tilde{\cal L}_t$. Now we have
$$\det({\cal E}_t\otimes{\cal L}_t)=\det({\cal E}_0),$$
in particular, ${\cal E}_t\otimes{\cal L}_t$ is a Ricci-flat deformation. Hence we obtain 
$${\cal E}_t\cong{\cal E}_0\otimes{\cal L}_t^\vee.$$

By imposing $H^2({\cal End}({\cal E}_0))=H^2({\cal O}_X)$ we ensure that every $\zeta\in H^1({\cal End}({\cal E}_0))$ can be integrated to a small
deformation (see \cite[Cor. 5.7]{b95}). So
$$\zeta=\frac{d}{dt}|_{t=0}{\cal L_t}^\vee\in H^1({\cal O}_X)$$
for some deformation ${\cal L}_t\in Pic^0(X)$ of ${\cal L}_0={\cal O}_X$. 
\end{proof}

\begin{prop}\label{hopfcoh}Let $X$ be a compact manifold allowing for a Ricci-rigid vector bundle and ${\cal E}_0$ a holomorphic vector bundle on $X$. 
If $H^1(ad({\cal E}_0))=0$, then ${\cal E}_0$ is Ricci-rigid.\end{prop}

\begin{proof}If $H^1({\cal End}({\cal E}_0))=H^1({\cal O}_X)$, and ${\cal E}_t$ is a Ricci-flat deformation, then by Grauert's semi-continuity theorem
we obtain also $H^1({\cal End}({\cal E}_t))=H^1({\cal O}_X)$ for small $t$. Shifting the centre of the deformation to $t$, we obtain a family $h(t)\in H^1({\cal End}({\cal E}_t))$.
We know now that $\tr:H^1({\cal End}({\cal E}_t))\seq H^1({\cal O}_X)$ is an isomorphism. By the Ricci-flatness of the deformation we have due to Proposition \ref{v=0} that
$\tr h(t)=0$. Hence $h(t)=0$ for all $t$ and so the deformation is trivial.
\end{proof}

\section{Examples}\label{examples}

% \subsection{K\"ahler manifolds}
% 
% \begin{cor}\label{kaehler}Every compact K\"ahler manifold is a ``stable curvature'' manifold. 
% \end{cor}
% 
% \begin{proof}Let $X$ be a compact K\"ahler manifold, $i:{\bb R}\seq{\cal O}_X$ be the inclusion and
% $j_p:{\cal A}^p_X\seq{\cal A}^{0,p}_X$ be the map mapping a real-valued $p$-form $\omega$ to its
% complex-valued $(0,p)$-part $\omega^{(0,p)}$ for $p\ge 1$ and the natural inclusion for $p=0$. Since the maps
% $$\xymatrix{0\ar[r]&{\bb R}\ar[r]\ar[d]_{i}&{\cal A}^0_X\ar[r]^d\ar[d]_{j_0}&{\cal A}^1_X\ar[r]^d\ar[d]_{j_1}&{\cal A}^2_X\ar[r]^d\ar[d]_{j_2}&\dots\\
% 0\ar[r]&{\cal O}_X\ar[r]&{\cal A}^{0,0}_X\ar[r]^{\overline\partial}&{\cal A}^{0,1}_X\ar[r]^{\overline\partial}&{\cal A}^{0,2}_X\ar[r]^{\overline\partial}&\dots}$$
% give a cochain map between these two acyclic resolutions, general theory (e.g. \cite[1,IV,4.4]{hs}) yields, 
% that $j_*=i_*:H^p(X,{\bb R})\seq H^p(X,{\cal O}_X)$
% on all levels $p$.
% 
% Any harmonic $(0,1)$-form $\eta$ yields by conjugation a harmonic $(1,0)$-form $\overline\eta$ such that $\omega:=\eta+\overline\eta$ is a
% real one-form with $i_*[\omega]=[\eta]$, so in the K\"ahler case $i_*$ is surjective.
% \end{proof}

Since all K\"ahler manifolds are stable curvature manifolds, Ricci-rigidity coincides with rigidity in these cases. Therefore we concentrate on non-K\"ahler examples.

\subsection{Hopf manifolds}

\begin{ex}{}\label{exhopf} Let $X$ be the Hopf manifold defined by the automorphism $\phi(z)=2z$ on ${\bb C}^n\setminus\{0\}$. Then there is a natural
smooth elliptic fibration $\pi:X\seq{\bb P}^{n-1}$. For any bundle we have $R^1\pi_*\pi^*{\cal E}_0={\cal E}_0^\vee$. 
Let ${\cal E}_0$ be a a simple bundle on ${\bb P}^{n-1}$ with
$H^1({\cal End}({\cal E}_0))=H^2({\cal End}({\cal E}_0))=0$. Then the Leray spectral sequence 
implies that $H^1({\cal End}(\pi^*{\cal E}_0))={\bb C}$ and hence $\pi^*{\cal E}_0$
is Ricci-rigid. For instance, $T_{{\bb P}^{n-1}}$ satisfies these conditions for $n\ge 2$. (For $n=2$ exactly the line bundles
satisfy the requirements.)
\end{ex}

Moreover, it is known that $p^*T_{{\bb P}^n}$ is not trivial for $n\ge 2$, if $p:{\bb C}^{n+1}\setminus\{0\}\seq{\bb P}^n$ denotes the natural projection. 
The following results will recover this and give a connection between Ricci-rigid bundles on some Hopf manifolds and non-trivial
vector bundles on ${\bb C}^n\setminus\{0\}$. 

\begin{prop}\label{hopf}Let $X$ be the Hopf manifold given by the quotient of ${\bb C}^n\setminus\{0\}$ by the automorphism group generated by 
$\varphi(z_1,\dots,,z_n)=(\alpha_1z_1,\dots,\alpha_nz_n)$, $|\alpha_i|>1$ and 
$u:{\bb C}^n\setminus\{0\}\seq X$ the projection. 
If ${\cal E}_0$ is a Ricci-rigid vector bundle on $X$ with $\rk{\cal E}_0>1$, then $u^*{\cal E}_0$ is not trivial.
\end{prop}

\begin{proof}
There is a multiplicative degree $\deg_\varphi:{\bb C}[z_1,\dots,z_n]\seq{\bb C}$ via
$$\deg_\varphi(z_i):=\alpha_i.$$

If $u^*{\cal E}_0$ is trivial, usual techniques allow us to identify ${\cal E}_0$ of rank $r$ with an equivalence class of holomorphic maps
$L:{\bb C^n}\setminus\{0\}\seq Gl(r,{\bb C})$ where
$$L\cong \tilde L:\iff\,\,\exists T\in{\cal O}({\bb C}^n\setminus\{0\}, Gl(r,{\bb C}))\mbox{ such that }
\tilde L=T\circ\varphi\cdot L\cdot T^{-1}.$$
Choosing $T$ carefully we can achieve a normal form of $L$ consisting of blocks $L_\nu, \nu\in{\bb C}^*$ in upper triangle form (cf. \cite{ma92} for a very similar normal form) with the property
$$(L_\nu)_{kk}=\prod_{j=1}^n\alpha_j^{i_{jk}}\nu$$ 
for $i_{jk}\ge 0, i_{j1}=0$ and 
$$(L_\nu)_{kl}\in{\bb C}[z_1,\dots,z_n]$$
homogeneous with 
$$\deg_\varphi(L_\nu)_{kl}=\prod_{j=1}^n\alpha_j^{i_{jl}-i_{jk}}.$$
Now $L(t)_{kl}:=L_{kl}\mbox{ for }(k,l)\notin\{(1,1),(2,2)\}$ and
$$L(t)_{11}:=\exp(t)L_{11}, L(t)_{22}:=\exp(-t)L_{22}$$
defines a non-trivial Ricci-flat small deformation of ${\cal E}_0$. 
\end{proof}

Combining (\ref{hopf}) and (\ref{exhopf}) we obtain immediately

\begin{cor}If $n>1, p:{\bb C}^{n+1}\setminus\{0\}\seq{\bb P}^n$ is the natural projection and ${\cal E}_0$ a simple vector bundle on ${\bb P}^n$ satisfying
$\rk{\cal E}_0>1$ and $H^1({\cal End}({\cal E}_0))=H^2({\cal End}({\cal E}_0))=0$, then $p^*{\cal E}_0$ is not trivial.
\end{cor}

\subsection{Surfaces}

\begin{thm}\label{nonfilt}Let $X$ be a compact complex surface with $b_1(X)\not=0$. 
If a rank 2 vector bundle ${\cal E}_0$ on $X$ is Ricci-rigid, then every filtration
$$0\seq{\cal F}\seq{\cal E}_0\seq{\cal G}\otimes I_Z\seq 0$$
with $Z$ a (non-reduced) point set and line bundles ${\cal F},{\cal G}$ satisfies 
$H^1(X,{\cal F}\otimes{\cal G}^\vee)\not=0.$
\end{thm}

\begin{proof}
We may assume that $X$ is not K\"ahler; otherwise there are no Ricci-rigid vector bundles. 
Assume we have an extension
$$0\seq{\cal F}\seq{\cal E}_0\seq{\cal G}\otimes I_Z\seq 0$$
with line bundles ${\cal F}, {\cal G}$ and a finite set of points $Z$ such that $H^1(X,{\cal F}\otimes{\cal G}^\vee)=0$.

This extension yields an element $\zeta\in\Ext^1(I_z,{\cal F}\otimes{\cal G}^\vee)$. A part of the long exact $Ext$-sequence
of the ideal sheaf sequence looks like
$$0=H^1(X,{\cal F}\otimes{\cal G}^\vee)
%\cong&\Ext^1({\cal O}_X,{\cal F}\otimes{\cal G}^\vee)
\seq\Ext^1(I_Z,{\cal F}\otimes{\cal G}^\vee)\seq\Ext^2({\cal O}_Z,{\cal F}\otimes{\cal G}^\vee)\seq H^2(X,{\cal F}\otimes{\cal G}^\vee).$$
So $\zeta$ may be identified with an element
$\xi\in\Ext^2({\cal O}_Z,{\cal G}^\vee\otimes{\cal F})$ mapping to $0\in H^2(X,{\cal G}^\vee\otimes{\cal F})$
and generating ${\cal Ext}^2({\cal O}_Z,{\cal F}\otimes{\cal G}^\vee)\cong{\cal O}_Z$. Now let
${\cal L}_t$ be a non-constant holomorphic family of line bundles with ${\cal L}_0={\cal O}_X$; due to $b_1(X)\not=0$ such a family exists. We interpret
$$\xi\in\Ext^2({\cal O}_Z,{\cal G}^\vee\otimes{\cal F})\cong\Ext^2({\cal O}_Z\otimes{\cal L}_t^{-2},{\cal G}^\vee\otimes{\cal F})\cong\Ext^2({\cal O}_Z,{\cal G}^\vee\otimes{\cal F}\otimes{\cal L}_t^2),$$
again generating ${\cal Ext}^2({\cal O}_Z,{\cal F}\otimes{\cal G}^\vee\otimes{\cal L}_t^2)\cong{\cal O}_Z$. 
By Grauert's semicontinuity and Lemma \ref{lbnonk} we may assume 
$$H^1(X,{\cal F}\otimes{\cal G}^\vee\otimes{\cal L}_t^2)=0$$
for all $t$ small enough and 
$$H^2(X,{\cal F}\otimes{\cal G}^\vee\otimes L_t^2)=0$$
for all $t\not=0$ small enough. So by the analogue $Ext$-sequence argument we obtain a unique $\zeta_t\in\Ext^1(I_Z,{\cal F}\otimes{\cal G}^\vee\otimes{\cal L}_t^2)$, i.e.
we obtain by the Serre construction (see e.g. \cite{b95}) an extension
$$0\seq{\cal F}\otimes{\cal L}_t\seq{\cal E}_t\seq{\cal G}\otimes{\cal L}_t^\vee\otimes I_Z\seq 0$$
with a vector bundle ${\cal E}_t$. 
The Ricci-flat deformation of ${\cal E}_0$ obtained in this way is not trivial, again due to Lemma \ref{lbnonk}: 
Of course, $H^0({\cal E}_0\otimes{\cal F}^\vee)\not=0$ and $H^0({\cal E}_t\otimes{\cal F}^\vee\otimes{\cal L}_t^\vee)\not=0$, hence $H^0({\cal E}_t\otimes{\cal F}^\vee)=0$ 
for general $t$, if the family ${\cal L}_t$ was chosen general enough. 
\end{proof}

\begin{lemma}\label{lbnonk}On a compact, non-K\"ahler surface $X$ for every line bundle ${\cal F}$ the set
$$\{{\cal L}\in Pic^0(X)|\,\,H^0({\cal F}\otimes{\cal L})=0\}$$
is open and not empty. The same holds true, if ${\cal F}$ is a filtrable rank $2$ vector bundle.
\end{lemma}

\begin{proof}First, let ${\cal F}$ be a line bundle. By Grauert's semicontinuity result it is enough to find some ${\cal L}\in Pic^0(X)$ with $H^0({\cal F}\otimes{\cal L})=0$.
So we assume that $H^0({\cal F}\otimes{\cal L})\not=0$ for all ${\cal L}\in Pic^0(X)$. Note that due to the non-K\"ahler property $c_1({\cal G})^2\le 0$ for every
line bundle ${\cal G}$ (cf. \cite[IV,Thm 6.2]{bhpv}) and $b_1(X)$ is odd (cf. \cite[IV,Thm 3.1]{bhpv}). The latter implies that $Pic^0(X)$ is not compact. 

In a first step we construct a fibration $\pi:X\seq C$. For $s\in H^0({\cal F})$ and $t\in H^0({\cal F}\otimes{\cal L})$ for some ${\cal L}\in Pic^0(X)$ we see
that due to
$$c_1({\cal F}).c_1({\cal F}\otimes{\cal L})=c_1({\cal F})^2\le 0$$
the schemes $V(s), V(t)$ do not intersect or have at least one common component. So
$$D:=\bigcap_{{\cal L}\in Pic^0(X)}V(H^0({\cal F}\otimes{\cal L}))$$
is empty or an effective divisor. The line bundle
$${\cal G}:={\cal F}\otimes{\cal O}_X(-D)$$
again satisfies $H^0({\cal G}\otimes{\cal L})\not=0$ for all ${\cal L}\in Pic^0(X)$ and, moreover, we can find ${\cal L}\in Pic^0(X)$ and sections
$s_0\in H^0({\cal G}), s_1\in H^0({\cal G}\otimes{\cal L})$ and $s_{-1}\in H^0({\cal G}\otimes{\cal L}^\vee)$ such that
$$V(s_0)\cap V(s_1)=V(s_0)\cap V(s_{-1})=\emptyset.$$
Hence, 
$$\pi':X\seq{\bb P}^1, x\mapsto [s_0^2:s_1s_{-1}]$$
is a well-defined map. Applying Stein factorization we obtain a fibration
$$\pi:X\seq B$$
onto a smooth curve $B$ with connected fibres. 

Note that all irreducible curves $C\subset X$ are contained in the fibres, since otherwise we have a curve $C$ with $C.F>0$
for the general fibre $F$ of $\pi$; but $(C+kF)^2=C^2+2kC.F>0$ for $k$ big enough. 

Let $C_1,\dots, C_n$ be the finitely many irreducible components of the reducible fibres of $\pi$. Now we have proved that for the compact manifold $M:=B\cup\{C_1,\dots, C_n\}$
and $P:={\cal F}\otimes Pic^0(X)$, the connected component of $Pic(X)$ containing ${\cal F}$, we have a map
$$\psi:\bigcup_{k\in{\bb N}_0}M^k\seq Pic(X)$$
induced by $\psi(b):={\cal O}_X(F_b)$ for $b\in B$ and $\psi(C_i):={\cal O}_X(C_i)$ satisfying 
$$P\subset\image\psi.$$

On the other hand $\psi(M^k)$ is compact and hence of codimension at least $1$ in $Pic(X)$. So $\image\psi$ cannot contain any connected component of $Pic(X)$
and we have reached a contradiction.

If ${\cal F}$ is a filtrable rank $2$ bundle, we take a filtration
$$0\seq{\cal F}'\seq{\cal F}\seq{\cal F}''\otimes I_Z\seq 0$$
and apply the lemma to the line bundles ${\cal F}'$ and ${\cal F}''$.
\end{proof}

\subsubsection{Ricci-rigid rank $2$ bundles on Inoue surfaces}

We treat Inoue surfaces before Hopf surfaces, since the description of Ricci-rigid vector bundles is somewhat easier due to the lack of curves.
The geometry of Inoue surfaces implies immediately that all line bundles are Ricci-rigid. So we tackle rank $2$ bundles.

\begin{cor}On an Inoue surface $X$ every Ricci-rigid rank $2$ vector bundle ${\cal E}_0$ is either non-filtrable or there
is a line bundle ${\cal F}$ such that ${\cal E}_0\otimes {\cal F}$ is isomorphic
to the unique vector bundle obtained as a non-split extension 
$$0\seq K_X\seq {\cal E}_0\otimes{\cal F}\seq{\cal O}_X\seq 0.$$ 
The latter is Ricci-rigid, indeed.
\end{cor}

\begin{proof}
By tensoring with line bundles we may assume that ${\cal E}_0$ allows for a filtration
$$0\seq{\cal L}\seq{\cal E}_0\seq I_Z\seq 0.$$
By Theorem \ref{nonfilt} we obtain $H^1(X,{\cal L})\not=0$. On an Inoue surface the intersection form of line
bundles is trivial, hence $\chi({\cal L})=\chi({\cal O}_X)=0$; furthermore 
$$H^0({\cal F})\not=0\iff {\cal F}={\cal O}_X$$
for any line bundle, since there are no curves on $X$. So the condition $H^1(X,{\cal L})\not=0$ translates
into ${\cal L}\in\{K_X,{\cal O}_X\}$.

First we treat the case ${\cal L}=K_X$ and $Z\not=\emptyset$. The filtration dualises to
$$0\seq{\cal O}_X\seq{\cal E}_0^\vee\seq K_X^\vee\otimes I_Z\seq 0.$$
So $H^0(X,{\cal E}_0^\vee)={\bb C}$ and, looking at the dualised sequence again, also $H^0(X,{\cal E}_0^\vee\otimes I_Z)={\bb C}$.
If we tensorise this sequence by ${\cal E}_0$ we obtain
$$0\seq{\cal E}_0\seq{\cal End}({\cal E}_0)\seq{\cal E}_0^\vee\otimes I_Z\seq 0,$$
yielding for global sections $H^0(X,{\cal End}({\cal E}_0))={\bb C}$, so ${\cal E}_0$ is simple. Tensoring the filtration
and its dual with $K_X$ tells us $H^0(X,{\cal E}_0\otimes K_X)=0=H^0(X,{\cal E}_0^\vee\otimes K_X\otimes I_Z)$, so also
$$H^2(X,{\cal End}({\cal E}_0)=H^0(X,{\cal End}({\cal E}_0)\otimes K_Z)=0.$$
Now we employ
$$-4c_2({\cal E}_0)=\chi({\cal End}({\cal E}_0))=1-h^1(X,{\cal End}({\cal E}_0)).$$
Theorem \ref{rrcoh} implies that ${\cal E}_0$ is not Ricci-rigid.

If $Z=\emptyset$ the filtration is
$$0\seq K_X\seq{\cal E}_0\seq{\cal O}_X\seq 0.$$
If the sequence splits, ${\cal E}_0$ is clearly not Ricci-rigid. For the non-splittung extension we have
$H^0({\cal E}_0)=0$.
%otherwise we would have a global section $t\in H^0(X,{\cal E}_0)$ such that for the section
%$s\in H^0(X,{\cal E}_0\otimes K_X)=H^0(X,{\cal E}_0^\vee)$ defining the filtration holds $s(t)\equiv 1$; the map
%${\cal O}_X\seq{\cal E}_0, f\mapsto ft$ would give a splitting.   
The dual of the filtration implies $H^0({\cal E}_0^\vee)={\bb C}$, so like above we obtain
that ${\cal E}_0$ is simple and $H^2(X,{\cal End}({\cal E}_0))=0$. Riemann-Roch implies this time
$H^1(X,{\cal End}({\cal E}_0))={\bb C}$ and so ${\cal E}_0$ is Ricci-rigid, indeed. 

We are left with the case $L={\cal O}_X$. Now ${\cal E}_0\cong{\cal E}_0^\vee$ and with the same techniques as above we see
that $H^2(X,{\cal End}({\cal E}_0))=0$ and $h^0({\cal End}({\cal E}_0))\in\{1,2\}$, if $Z\not=\emptyset$. Riemann-Roch
yields again that ${\cal E}_0$ is not Ricci-rigid. If $Z=\emptyset$, the non-splitting extension
$$0\seq{\cal O}_X\seq{\cal E}_0\seq{\cal O}_X\seq 0$$
satisfies $H^0(X,{\cal E}_0)={\bb C}$. Hence $h^0(X,{\cal End}({\cal E}_0))\ge 2$.
%any section $s\in H^0(X,{\cal E}_0)$ gives an endomorphism $s\otimes s$ of rank $1$. %%% E=E dual! 
Moreover, $$H^1(X,{\cal E}_0)={\bb C}, H^2(X,{\cal E}_0)=0.$$ The sequence
$$0\seq{\cal E}_0\seq{\cal End}({\cal E}_0)\seq{\cal E}_0\seq 0$$
now implies $H^2(X,{\cal End}({\cal E}_0))=0, h^0(X,{\cal End}({\cal E}_0))=h^1(X,{\cal End}({\cal E}_0))=2$,
so ${\cal E}_0$ is not Ricci-rigid.
\end{proof}

\subsubsection{Ricci-rigid vector bundles on Hopf surfaces}

\begin{cor}Let $X$ be a Hopf surface given by diagonal $\varphi$. 
All line bundles on $X$ are Ricci-rigid. There are no Ricci-rigid vector bundles of rank $\ge 3$. There
are no good Ricci-rigid rank $2$ vector bundles.
\end{cor}

\begin{proof}By Proposition \ref{hopfcoh} the line bundles are Ricci-rigid.
 
Next we observe that a vector bundle ${\cal G}$ on ${\bb C^2}\setminus\{0\}$ is trivial if and only if it is extendable to 
${\bb C}^2$ (see \cite{s66}).
It is extendable as a vector bundle if and only if it is extendable as a coherent sheaf, but this is ensured by \cite[Cor. VII.4]{fg69}, if $\rk{\cal G}\ge 3$. Hence, by
Proposition \ref{hopf}, a vector bundle
${\cal E}_0$ on a Hopf surface given by diagonal $\varphi$ is never Ricci-rigid, if $\rk{\cal E}_0>2$. So it remains to deal with the ``critical case'' $\rk{\cal E}_0=2$.   

In case ${\cal E}_0$ is Ricci-rigid, good and $\rk{\cal E}_0=2$ we have $H^1(X,{\cal End}({\cal E}_0))={\bb C}$ and
deduce
$$-4c_2({\cal E}_0)=\chi({\cal End}({\cal E}_0))\ge 0.$$
If $c_2({\cal E}_0)<0$, then ${\cal E}_0$ is not filtrable, so it is simple;
%see B95, p.91 Def 4.9: If E is not simple, there is a non-invertible non-zero endomorphism; its kernel gives
%a filtration.
but then $c_2({\cal E}_0)=0$.
If $c_2({\cal E}_0)=0$, then ${\cal E}_0$ is simple and topologically trivial, contradicting \cite[Prop 3.1]{mo04}.
\end{proof}


\begin{thebibliography}{BrMo05}
\bibitem[AT03]{at03}{Aprodu, M., Toma, M.}, 'Une note sur les fibr\'es holomorphes non-filtrables. (French. English, French summary) [A note on non-filtrable holomorphic bundles]',
{\em C. R. Math. Acad. Sci. Paris} 336(7), 581--584 (2003) 
\bibitem[Ba02]{b02}{Ballico, E.}, 'Holomorphic vector bundles on ${\bb C}^2\setminus\{0\}$', {\em Isr. J. of Math.} 128, 197--204 (2002)
\bibitem[BP87]{bp87}{B\u anic\u a, C., Le Potier, J.}, 'Sur l'existence des fibrés vectoriels holomorphes sur les surfaces non-algébriques', {\em J. Reine Angew. Math.} 378, 1--31 
(1987) 
\bibitem[BHPV]{bhpv}{Barth, W., Hulek, K., Peters, C., van de Ven, A.}, 'Compact Complex Surfaces', Springer, 2004
\bibitem[Br95]{b95}{Br\^inz\u anescu, V.}, 'Holomorphic Vector Bundles over Compact Complex Surfaces', Springer, 1996
\bibitem[BrMo05]{bm05}{Br\^inz\u anescu, V., Moraru, R.}, 'Stable Bundles on Non-K\"ahler Elliptic Surfaces', {\em Commun. Math. Phys.} 254, 565--580 (2005)
\bibitem[FG69]{fg69}{Frisch, J.; Guenot, J.}, 'Prolongement de faisceaux analytiques coh\'erents', {\em Invent. Math.}  7,  321--343 (1969)
\bibitem[Ha]{ha}{Hartshorne, R.}, 'Algebraic Geometry', Springer, 1977
\bibitem[HS]{hs}{Hilton, P.J., Stammbach, U.}, 'A Course in Homological Algebra', Springer, 1971
\bibitem[Ho64]{ho64}{Horrocks, G.}, 'Vector bundles on the punctured spectrum of a local ring', {\em Proc. London Math. Soc. (3)} 14, 689--713 (1964) 
\bibitem[Mo04]{mo04}{Moraru, R.}, 'Stable bundles on Hopf manifolds', math.AG/0408439 (2004)
\bibitem[Ma92]{ma92}{Mall, D.}, 'On holomorphic vector bundles on Hopf manifolds with trivial pullback on the universal covering', {\em Math. Ann.} 294, 718--740 (1992) 
\bibitem[S66]{s66}{Serre, J.-P.}, 'Prolongement de faisceaux analytiques coh\'erents', {\em Ann. Inst. Fourier} 16, No.1, 363--374 (1966)
\bibitem[Se58]{se58}{Seshadri, C.S.}, 'Triviality of vector bundles over the affine space $K^2$', {\em Proc. Natl. Acad. Sci. USA } 44, 456--458 (1958)
\bibitem[Ve06]{v06}{Verbitsky, M.}, 'Holomorphic bundles on diagonal Hopf manifolds', {\em Izvestiya} 70(5), 867--882 (2006)
\end{thebibliography}
\end{document}